\documentclass[3p]{elsarticle}

\usepackage{amsmath}
\usepackage{amsthm}
\usepackage{amssymb}
\usepackage{bm}

\newdefinition{rmk}{Remark}
\newproof{pf}{Proof}

\graphicspath{{figs/}}

\journal{arXiv}

\begin{document}

\begin{frontmatter}

\title{Approximate representation of the solutions of fractional elliptical BVP through the solution of parabolic IVP\tnoteref{label1}}
\tnotetext[label1]{The publication has been prepared with the support of the \emph{RUDN University Program 5-100}.}

\author{P.N. Vabishchevich\corref{cor1}\fnref{lab1,lab2}}
\ead{vabishchevich@gmail.com}
\cortext[cor1]{Correspondibg author.}

\address[lab1]{Nuclear Safety Institute, Russian Academy of Sciences,
              52, B. Tulskaya, 115191 Moscow, Russia}

\address[lab2]{Peoples' Friendship University of Russia (RUDN University), 6 Miklukho-Maklaya St, 117198 Moscow, Russia}

\begin{abstract}
Boundary value problem for a fractional power of an elliptic operator is considered. 
An integral representation by means of a standard solution problem for parabolic equations is used to solve such problems. Quadrature generalized Gauss-Laguerre formulas are constructed. 
We examine the effect of key parameters on the accuracy of the approximate solution: the number of nodes of the quadrature and fractional power of the operator. 
Computational experiments were performed to model two-dimensional problem with a fractional power of an elliptic operator.
\end{abstract}


\end{frontmatter}

\section{Introduction} 
\label{sec:1}

Various nonlocal process models \cite{baleanu2012fractional,uchaikin}, which are based on the use of fractional power of an elliptic operator \cite{Pozrikidis16}, are currently actively discussed. 
In the simplest case, the problem for the fractional Laplacian is formulated as follows.
For a bounded domain $\Omega$ on the set of functions $v(\bm x) = 0, \ \bm x \in \partial \Omega$ we define the linear operator $\mathcal{A}: \  \mathcal{A} v = - \triangle v, \ \bm x \in \Omega$.
We seek the solution of the problem for the equation with the fractional power elliptic operator
$\mathcal{A}^{\alpha } v = f, \ 0 < \alpha < 1$ for a given $f(\bm x), \ \bm x \in \Omega$.

Various computational approaches for approximate solution of the problems with fractional power of an operator are used. As applied to problems of linear algebra, they are discussed in the work \cite{higham2008functions}.
For the problems with the fractional power of the elliptic operator under the consideration, special attention should be paid to the methods whose computational implementation is based on a solution of standard boundary value problems \cite{bonito2017}.

For a functional approximation of fractional powers, best uniform rational approximation \cite{stahl2003best} can be applied.
This approach for approximate solving of multidimensional problems of fractional diffusion is implemented in the works
\cite{harizanov2018optimal,HarizanovArxive2019}. 
The fractional power of the operator is approximated by the sum of standard operators. 
Similar approximations can be obtained using the appropriate integral representation of fractional powers of the operator and various quadrature formulas.
Classical integral Balakrishnan representations \cite{balakrishnan1960fractional,birman1987spectral} are the basis of the works 
\cite{bonito2015numerical,AcetoNovat,Aceto2019}.
Other integral representations are used in our paper \cite{vab2019-10838}.

A solution to a problem with a fractional power of an elliptic operator can be obtained from a solution to a problem of a larger dimension.
In \cite{Caffarelli} it is shown that the solution of the fractional Laplacian problem can be obtained as a solution to the elliptic problem on the semi-infinite cylinder domain. Computational algorithms based on such a relationship are discussed, for example, in
\cite{nochetto2015pde}.

For solving fractional power elliptic problems, we have proposed  \cite{vabishchevich2014numerical} a numerical algorithm based on a transition to a pseudo-parabolic equation (Cauchy problem method).  The computational algorithm is simple for practical use, robust, and applicable to solving a wide class of problems. To increase the accuracy of the approximate solution, two- and three-level schemes of higher order approximation on regular and irregular grids (see, for example,
\cite{duan2018numerical,ciegisvab2019-00201}) are used.

Boundary value problems for the fractional power of an elliptic operator are often considered \cite{stinga2019user} in terms of the methods of semigroups.
In the works \cite{cusimano2018discretizations, cusimano2018numerical}
, a singular integral expression is used to construct an approximate solution by solving the standard Cauchy problem for the parabolic equation. Special quadrature formulas that take into account the singularity and finite element approximations in space are used.

In the present paper, we also use the representation of the solution of the boundary value problem for the fractional power of a positive definite elliptic operator by solving IVP for the parabolic equation. The amplitude dependence of the solution of the auxiliary parabolic problem is highlighted. Approximations of the desired solution of fractional elliptical BVP is based on generalized Gauss-Laguerre formula \cite{davis2007methods}, which most fully takes into account the singularity of the integral representation.

The paper is organized as follows. Section 2 describes a problem for an equation with a fractional power of an elliptic operator of second order. The problem is considered in a rectangle using a standard finite-difference approximation on a uniform grid. Singular integral expression through the solution of a standard Cauchy problem for the parabolic equation is discussed in Section 3. In Section 4, we construct the corresponding quadrature formulas. Numerical experiments for a model two-dimensional problem with the fractional power of the Laplace operator are given in Section 5. The results of the work are summarized in Section 6.

\section{Problem formulation}\label{sec:2}

In a bounded polygonal domain $\Omega \subset R^2$ with the Lipschitz continuous boundary $\partial\Omega$,
we search the solution of a  problem with a fractional power of an elliptic operator.
Here, we use the definition of a fractional power of an elliptic operator that relies on the spectral 
theory \cite{birman1987spectral}.
The elliptic operator is introduced as
\begin{equation}\label{1}
  \mathcal{A}  v = - {\rm div}  ( a({\bm x}) {\rm grad} \, v) + c({\bm x}) v
\end{equation} 
with coefficients $0 < a_1 \leq a({\bm x}) \leq a_2$, $c({\bm x}) \geq 0$.
The operator $\mathcal{A}$ is defined on the set of functions $v({\bm x})$ that satisfy the following conditions on the boundary $\partial\Omega$:
\begin{equation}\label{2}
  v ({\bm x}) = 0,
  \quad {\bm x} \in \partial \Omega .
\end{equation} 

In the Hilbert space $\mathcal{H} = L_2(\Omega)$, we define 
the scalar product and norm in a standard way:
\[
  (v,w) = \int_{\Omega} v({\bm x}) w({\bm x}) d{\bm x},
  \quad \|v\| = (v,v)^{1/2} .
\] 
In the spectral problem
\[
 \mathcal{A}  \varphi_k = \lambda_k \varphi_k, 
 \quad \bm x \in \Omega , 
\] 
\[
  \varphi_k = 0,
  \quad {\bm x} \in \partial \Omega , 
\] 
we have 
\[
 0 < \lambda_1 \leq \lambda_2 \leq ... ,
\] 
and the eigenfunctions  $ \varphi_k, \ \|\varphi_k\| = 1, \ k = 1,2, ...  $ form a basis in $L_2(\Omega)$. Therefore, 
\[
 v = \sum_{k=1}^{\infty} (v,\varphi_k) \varphi_k .
\] 
Let the operator $\mathcal{A}$ be defined in the following domain:
\[
 D(\mathcal{A} ) = \Big \{ v \ | \ v(\bm x) \in L_2(\Omega), \ \sum_{k=0}^{\infty} | (v,\varphi_k) |^2 \lambda_k < \infty \Big \} .
\] 
Under these conditions  $\mathcal{A} : L_2(\Omega) \rightarrow L_2(\Omega)$ and
the operator $\mathcal{A}$ is self-adjoint and positive definite: 
\begin{equation}\label{3}
  \mathcal{A}  = \mathcal{A} ^* \geq \nu   \mathcal{I} ,
  \quad \nu   > 0 ,    
\end{equation} 
where $\mathcal{I}$ is the identity operator in $\mathcal{H}$.
In applications, the value of $\lambda_1$ is unknown (the spectral problem must be solved).
Therefore, we suppose that $\nu  \leq \lambda_1$ in (\ref{3}).
Let us assume for the fractional power of the  operator $\mathcal{A}$:
\[
 \mathcal{A} ^\alpha v =  \sum_{k=0}^{\infty} (v,\varphi_k) \lambda_k^\alpha  \varphi_k .
\] 
The solution $v(\bm x)$ satisfies the equation
\begin{equation}\label{4}
  \mathcal{A}^\alpha v = f 
\end{equation} 
under the restriction $0 < \alpha < 1$. 

We consider the simplest case where the computational domain $\Omega$ is a rectangle
\[
 \Omega = \{ \bm x  \ | \ \bm x = (x_1,x_2), \ 0 < x_n < l_n, \ n = 1,2 \} .
\]
To solve the problem approximately (\ref{4}), we introduce a uniform grid in the domain $\Omega$
\[
\overline{\omega}  = \{ \bm{x} \ | \ \bm{x} =\left(x_1, x_2\right), \quad x_n =
i_n h_n, \quad i_n = 0,1,...,N_n,
\quad N_n h_n = l_n, \ n = 1,2 \} ,
\]
where $\overline{\omega} = \omega \cup \partial \omega$ and
$\omega$ is the set of interior nodes, whereas $\partial \omega$ is the set of boundary nodes of the grid.
For the grid functions $u(\bm x)$ such as $u(\bm x) = 0, \ \bm x \notin \omega$, we define the Hilbert space
$H=L_2\left(\omega\right)$, where the scalar product and the norm are specified as follows:
\[
\left(u, w\right) =  \sum_{\bm x \in  \omega} u\left(\bm{x}\right)
w\left(\bm{x}\right) h_1 h_2,  \quad 
\| y \| =  \left(y, y\right)^{1/2}.
\]

For $u(\bm x) = 0, \ \bm x \notin \omega$, the grid operator $A$  can be written as
\[
  \begin{split}
  A u = & -
  \frac{1}{h_1^2} a(x_1+0.5h_1,x_2) (u(x_1+h_1,h_2) - u(\bm{x})) \\ 
  & + \frac{1}{h_1^2} a(x_1-0.5h_1,x_2) (u(\bm{x}) - u(x_1-h_1,h_2)) \\
  & - \frac{1}{h_2^2} a(x_1,x_2+0.5h_2) (u(x_1,x_2+h_2) - u(\bm{x})) \\ 
  & + \frac{1}{h_2^2} a(x_1, x_2-0.5h_2) (u(\bm{x}) - u(x_1,x_2-h_2)) + c(\bm x) u(\bm x), 
  \quad \bm{x} \in \omega . 
 \end{split} 
\] 
For the above grid operators (see \cite{Samarskii1989,SamarskiiNikolaev1978}), we have 
\begin{equation}\label{5}
 A = A^* \geq \delta I,
 \quad \delta > 0 ,
\end{equation} 
where $I$ is the grid identity operator in $H$. 
For problems with sufficiently smooth coefficients and the right-hand side, 
it approximates the differential operator with the truncation error 
$\mathcal{O} \left(|h|^2\right)$, $|h|^2 = h_1^2+h_2^2$. 

To handle the fractional power of the grid operator $A$, let us consider the eigenvalue problem
\begin{equation}\label{6}
 A \psi_k = \mu_k \psi_k .  
\end{equation} 
We have
\[
 0 < \delta = \mu_1 \leq \mu_2 \leq ... \leq \mu_K,
 \quad K = (N_1-1)(N_2-1) , 
\] 
where eigenfunctions $\psi_k, \ \|\psi_k\| = 1, \ k = 1,2, ..., K,$ form a basis in $H$. Therefore
\[
 u = \sum_{k= 1}^{K}(u, \psi_k) \psi_k . 
\]
For the fractional power of the operator $A$, we have
\[
 A^\alpha u = \sum_{k= 1}^{K}(u, \psi_k) \mu_k^\alpha \psi_k .
\] 
Using the above approximations, from (\ref{4}), we arrive at the discrete problem
\begin{equation}\label{7}
 u = A^{-\alpha} b ,
\end{equation} 
where, for example, $b(\bm x) = f(\bm x) , \ \bm x \in \omega$.
This problem is considered for (\ref{5}), moreover, $\delta \leq \mu_1$.
 
\section{Representation by solving parabolic IVP}\label{sec:3}

An approximate solution to the boundary value problem for the fractional degree of an elliptic operator (\ref{7})
often based on the use of one or another integral representation.
For example, in the works mentioned above \cite{bonito2015numerical, AcetoNovat}
the Balakrishnan formula \cite{balakrishnan1960fractional, birman1987spectral} is used when
\begin{equation}\label{8}
 A^{-\alpha} = \frac{\sin(\pi \alpha)}{\pi} \int_{0}^{\infty} \theta^{-\alpha} (A + \theta  I)^{-1} d \theta ,
 \quad 0 < \alpha < 1 .
\end{equation} 
The approximation of $A^{-\alpha}$ is based on the use of one or another quadrature formulas for
the right-hand side of (\ref{8}). 

Here, the possibilities of approximating both the fractional power of the operator $A^{-\alpha}$ and approximating the solution of the problem are realized.
\begin{enumerate}
 \item Rational function approximation  $A^{-\alpha}$:
\[
 A^{-\alpha} \approx \sum_{i = 1}^{m} \gamma_i (A + \theta_i I)^{-1}, 
\] 
where $\gamma_i> 0$ are the weights, and $\theta_i> 0, \ i = 1, 2, \ldots, m$ are the nodes of the quadrature formula for (\ref{8}).
 \item An approximate solution to the problem (\ref{4}):
\[
 u \approx \sum_{i = 1}^{m} \gamma_i u_i,
 \quad (A + \theta_i I) u_i = b,
 \quad  i = 1, 2, \ldots, m .
\]  
\end{enumerate} 
The inverse operator of the fractional power elliptic problem 
is treated as a sum of inverse operators of classical elliptic operators
and the required approximate solution of the problem is taken as the sum of the solutions of standard elliptical BVP.

The integral representation of the fractional power of the elliptic operator, which is fundamentally different from (\ref{8})
used in works \cite{cusimano2018discretizations,cusimano2018numerical}.
It is based on the method 
of semigroups, когда $A^{-\alpha}$ is defined through 
an integral formulation as follows
\begin{equation}\label{9}
  A^{-\alpha} = \frac{1}{\Gamma(\alpha)} 
	\int_0^\infty \theta^{\alpha-1} e^{{- \theta  A}}   d \theta ,
  \quad 0 < \alpha < 1 ,
\end{equation}
where $\Gamma(\alpha)$ is gamma function. 
In this case, the solution to the problem (\ref{7}) is
\begin{equation}\label{10}
  u = \frac{1}{\Gamma(\alpha)} 
	\int_0^\infty \theta^{\alpha-1} e^{{- \theta  A}} b \, d \theta .
\end{equation}
Function $w(t) = e^{{- t A}} b$ is a solution 
of the following Cauchy problem
\begin{equation}\label{11}
 \frac{d w}{d t} + A w = 0,
 \quad 0 < t < \infty ,
\end{equation} 
\begin{equation}\label{12}
 w(0) = b .
\end{equation} 
By the representation (\ref{10}) we have
\begin{equation}\label{13}
  u = \frac{1}{\Gamma(\alpha)} 
	\int_0^\infty \theta^{\alpha-1} w(\theta) d \theta .
\end{equation}

For (\ref{1}), (\ref{2}) the problem (\ref{11}), (\ref{12}) is the standard
Cauchy problem for a second-order parabolic equation.
Based on (\ref{13}), the fractional elliptical BVP solution is represented by the IVP solution for the parabolic equation. Nonlocal stationary problem is associated with solving a non-stationary local problem and therefore, the considered approach is methodologically close to Cauchy problem method, which is proposed in \cite{vabishchevich2014numerical}.

We note the most important features of the approximation of the fractional power of the operator $A^{-\alpha}$
and the problem solutions (\ref{7}) based on the integral representation (\ref{9}).
\begin{enumerate}
 \item With quadrature formulas for (\ref{9}) we come to the approximation 
$A^{-\alpha}$ with operator exponentials:
 \[
  A^{-\alpha} \approx \sum_{i = 1}^{m} \gamma_i e^{- \theta_i A} . 
 \] 
Here, as before, $\gamma_i> 0$, $\theta_i> 0, \ i = 1, 2, \ldots, m$ are the weights and nodes of the corresponding quadrature formula.
 \item For an approximate solution to the problem (\ref{7}) from (\ref{13}) we have
\begin{equation}\label{14}
  u \approx u_m = \sum_{i = 1}^{m} \gamma_i w(\theta_i),
\end{equation} 
where $w(t)$ is the solution to the problem (\ref{11}), (\ref{12}).
\end{enumerate}  

Thus, the computational algorithm for solving the problem (\ref{7})
can be built on a numerical solution of the Cauchy problem (\ref{11}), (\ref{12}), according to which we have the approximate solution with the selected quadrature formula (see (\ref{14})).

The integral representation (\ref{13}) is singular, that is, the integrand has a singularity for $t \rightarrow 0$ for $0 <\alpha <1$.
For this reason, we can not rely on the accuracy of quadrature formulas, especially when $\alpha \rightarrow 0 $.
We extend the class of integral representations to eliminate the singularities of the integrand.

It is possible to transform (\ref{13}) taking into account (\ref{11}), (\ref{12}).
Based on the solutions of the Cauchy problem (\ref{11}), (\ref{12}) we have
\[
 \alpha \theta^{\alpha-1} w(\theta) d \theta =  d(\theta^{\alpha} w(\theta)) - \theta^{\alpha} \frac{d w}{d \theta} d \theta 
 =  d(\theta^{\alpha} w(\theta)) + \theta^{\alpha} A w .
\]  
Integrating by parts allows writing (\ref{13}) as
\begin{equation}\label{15}
  u = \frac{1}{\alpha \Gamma(\alpha)} 
	\int_0^\infty \theta^{\alpha} v(\theta) d \theta ,
\end{equation}
where $v = A w$. From (\ref{11}), (\ref{12}) we have the Cauchy problem for $v$:
\begin{equation}\label{16}
 \frac{d v}{d t} + A v = 0,
 \quad 0 < t < \infty ,
\end{equation} 
\begin{equation}\label{17}
 v(0) = A b .
\end{equation} 

We arrive at a similar result based on the representation (\ref{9}).
Replacing $\alpha$ by $\alpha + p$, where $p$ has a non-negative integer, we get
\begin{equation}\label{18}
  A^{-\alpha} = \frac{A^p}{\Gamma(\alpha+p)} 
	\int_0^\infty \theta^{\alpha+p-1} e^{{- \theta  A}}   d \theta ,
  \quad 0 < \alpha < 1 ,
\end{equation}
For the problem solution with a fractional power of the operator (\ref{5}), it follows from (\ref{18}) that
\begin{equation}\label{19}
  u = \frac{1}{\Gamma(\alpha+p)} 
	\int_0^\infty \theta^{\alpha+p-1} e^{{- \theta  A}} A^p b \, d \theta .
\end{equation}
We define function $w(t) = e^{{- t A}} A^p b$ as a solution of the equation (\ref{11}),
supplemented by the initial condition
\begin{equation}\label{20}
 w(0) = c,
 \quad c = A^p b .
\end{equation} 
To solve the problem (\ref{5}) from (\ref{19}) we obtain the representation
\begin{equation}\label{21}
  u = \frac{1}{\Gamma(\alpha+p)} 
	\int_0^\infty \theta^{\alpha+p-1} w(\theta) d \theta .
\end{equation}
For $p = 1$, the integral representation (\ref{21}) through the solution of the problem (\ref{11}), (\ref{20})
corresponds to (\ref{16})--(\ref{18}).

\section{Numerical integration}\label{sec:4} 

We have two sources of errors in the approximate representation of the solution of the boundary value problem
for the fractional power of a discrete elliptic operator through the solution of a parabolic problem.
The first is related to the approximate calculation of the integral (\ref{21}), the second is connected with a numerical solution of the Cauchy problem (\ref{11}), (\ref{20}). 
At this stage of our study, we assume that the problem (\ref{11}), (\ref{20}) is solved accurately and we will focus on the errors of approximate integration.

When using quadrature formulas (\ref{14}) for an approximate solution of the problem (\ref{7})
it is necessary to take into account the possible singularity of the integral representation: integration on the half-line and the features of the integrand.
In the works of Cusimano et al. \cite{cusimano2018discretizations, cusimano2018numerical}
integration (\ref{13}) is held at a uniform partitioning of the final interval and the weights of the quadrature formula are associated with the function $\theta^{\alpha -1}$.
Taking into account (\ref{6}) to solve the problem (\ref{11}), (\ref{20}) we have the representation
\[
 w(t) = \sum_{k= 1}^{K}(c, \psi_k) e^{-\mu_k t} \psi_k . 
\] 
From this perspective, from (\ref{21}) we obtain
\begin{equation}\label{22}
 u = \frac{1}{\Gamma(\alpha+p)} \sum_{k= 1}^{K}(c, \psi_k)  \psi_k \int_0^\infty \theta^{\alpha+p-1} e^{-\mu_k \theta} d \theta .
\end{equation} 
Introducing a new integration variable $\xi = \delta \theta$ from (\ref{22}) we get
\[
 u = \frac{1}{\delta^{\alpha+p} \Gamma(\alpha+p)} \sum_{k= 1}^{K}(c, \psi_k)  \psi_k \int_0^\infty
 \xi^{\alpha+p-1} e^{-\xi} e^{-\delta^{-1} (\mu_k - \delta)\xi } d \xi .
\] 
In view of this, it is necessary to construct quadrature formulas for the integrals
\begin{equation}\label{23}
 S(\alpha+p,\kappa) = \int_0^\infty \xi^{\alpha+p-1} e^{-\xi} \varrho (\xi; \kappa) d \xi ,
 \quad 0 < \alpha < 1 ,
 \quad p \geq 0 , 
\end{equation} 
in which the integrand is
\begin{equation}\label{24}
 \varrho (\xi; \kappa) = e^{- \kappa \xi}, 
 \quad \kappa \geq 0,  
\end{equation} 
Concerning (\ref{22}) we have
\[
 \kappa = \kappa_k,
 \quad \kappa_k = \frac{\mu_k }{\delta } - 1, 
 \quad k = 1,2, \ldots, K . 
\] 

When calculating the integrals (\ref{23}) with the weight function $\xi^{\alpha + p-1} e^{- \xi}, \ 0 <\alpha <1, \ p \geq 0$, we focus on the generalized Gauss-Laguerre quadrature formula \cite{davis2007methods}.
In this case
\begin{equation}\label{25}
 S(\alpha+p,\kappa) \approx S_m(\alpha+p,\kappa) = \sum_{i=1}^{m} \sigma_i \varrho (\xi_i; \kappa) ,
\end{equation} 
where the nodes of the quadrature formula $\xi_i$ are the zeros of the generalized Laguerre polynomial $L_m^{(\alpha + p -1)} (\xi)$,
and weights are 
\[
 \sigma_i = \frac{\Gamma(m+\alpha+p)}{m!} \xi_i \Big(L_{m+1}^{(\alpha+p -1)} (\xi_i) \Big)^2 ,
 \quad i = 1,2, \ldots, m . 
\] 

The accuracy of the approximate calculation of $S(\alpha + p, \kappa)$ is estimated at $\kappa \in [0, 10^5]$ by the value of the relative error
\[
 \varepsilon = \frac{1}{q} 
\max_{0 \leq \kappa \leq 10^5} |S(\alpha+p,\kappa) - S_m(\alpha+p,\kappa) | ,
 \quad q = \max_{0 \leq \kappa \leq 10^5} S(\alpha+p,\kappa) .
\] 
For convenience of the analysis of the calculated data, the dependence on $\alpha$ and $p$ is highlighted separately.

The accuracy of the generalized Gauss-Laguerre quadrature formula for calculating the integral $S(\alpha + p, \kappa)$
with a different number of nodes is given in Table 1.
The accuracy decreases with decreasing fractional power $\alpha$ and monotonically increases with increasing $\alpha + p$.
A fundamental increase in accuracy is achieved by choosing $p> 0$.

\begin{center}
\begin{table}[htp]
\label{t-1}
\caption{Approximation error $\varepsilon$ in the calculation $S(\alpha+p,\kappa)$}
\centering
\begin{tabular}{ccccccc}
\hline
   $m$   &    $p$   &   $\alpha = 0.1$         &     $\alpha = 0.25$   &    $\alpha = 0.5$   &     $\alpha = 0.75$      &    $\alpha = 0.9$   \\
\hline
        &      0   &  5.241730e-01     &  2.617133e-01     &    1.032809e-01    &  4.721032e-02    &  3.089846e-02  \\
        &      1   &  1.827771e-02     &  1.264074e-02     &    7.111007e-03    &  4.169156e-03    &  3.078424e-03  \\
   25   &      2   &  2.090947e-03     &  1.583269e-03     &    1.016776e-03    &  6.683762e-04    &  5.249533e-04  \\
        &      3   &  3.846239e-04     &  3.069550e-04     &    2.136824e-04    &  1.511406e-04    &  1.236642e-04  \\
        &      4   &  9.538633e-05     &  7.894328e-05     &    5.818673e-05    &  4.339988e-05    &  3.660028e-05  \\
\hline
        &      0   &  4.891202e-01     &  2.202112e-01     &    7.321202e-02    &  2.822765e-02    &  1.669060e-02  \\
        &      1   &  8.628584e-03     &  5.396535e-03     &    2.569714e-03    &  1.276760e-03    &  8.540722e-04  \\
   50   &      2   &  5.088380e-04     &  3.493823e-04     &    1.907784e-04    &  1.067459e-04    &  7.615334e-05  \\
        &      3   &  4.911145e-05     &  3.563271e-05     &    2.118067e-05    &  1.280562e-05    &  9.540530e-06  \\
        &      4   &  6.498686e-06     &  4.901612e-06     &    3.097237e-06    &  1.982550e-06    &  1.526013e-06  \\
\hline
        &      0   &  4.563886e-01     &  1.852328e-01     &    5.183319e-02    &  1.683123e-02    &  8.980301e-03  \\
        &      1   &  4.049587e-03     &  2.286550e-03     &    9.186559e-04    &  3.853388e-04    &  2.329379e-04  \\
  100   &      2   &  1.212817e-04     &  7.529116e-05     &    3.476968e-05    &  1.646274e-05    &  1.062789e-05  \\
        &      3   &  6.001109e-06     &  3.942100e-06     &    1.986278e-06    &  1.018559e-06    &  6.876088e-07  \\
        &      4   &  4.108546e-07     &  2.809418e-07     &    1.508100e-07    &  8.206209e-08    &  5.731183e-08  \\
\hline
\end{tabular}
\end{table}
\end{center}

\section{The numerical solution of fractional elliptical BVP}\label{sec:5} 

We will illustrate the possibilities of an approximate solution to the problem (\ref{1}), (\ref{2}), (\ref{4})
based on the solution of the Cauchy problem (\ref{11}), (\ref{20}) using the quadrature formula (\ref{25}).
The model problem for the Laplace operator is considered when
\[
 a(\bm x) = 1,
 \quad c(\bm x) = 0,
 \quad \bm x \in \Omega .  
\]   

For eigenfunctions and eigenvalues (\ref{6}) we have (see, e.g., \cite{SamarskiiNikolaev1978}):
\[
\begin{split}
 \psi_k (\bm x) & = \prod_{\beta  =1}^{2} \sqrt{\frac{2}{l_\beta } } \sin (k_\beta  \pi x_\beta ),
 \quad \bm x \in \omega , \\ 
 \mu_k & = \sum_{\beta  =1}^{2}  \frac{4}{h_\beta^2} \sin^2 \frac{k_\beta  \pi}{2 N_\beta } ,
 \quad k_\alpha  = 1,2,...,N_\alpha -1, \quad  \alpha  = 1, 2 .
\end{split} 
\] 
Direct calculations yield 
\[
 \mu_1 = \sum_{\beta  =1}^{2}  \frac{4}{h_\beta ^2} \sin^2 \frac{\pi}{2 N_\beta } 
 <  8 \Big (\frac{1}{l_1^2} + \frac{1}{l_2^2} \Big ) ,
 \quad \lambda_K = \sum_{\beta  =1}^{2}  \frac{4}{h_\beta ^2} \cos^2 \frac{\pi}{2 N_\beta } 
 <  4 \Big (\frac{1}{h_1^2} + \frac{1}{h_2^2} \Big ) .
\] 

The model problem (\ref{4}) is solved numerically with two right-hand sides: $f(\bm x) = f_\beta (\bm x), \ \beta = 1,2$, when
\[
 f_1 (\bm x) = x_1^2 (1-x_1) x_2^2(1-x_2),
 \quad  f_2 (\bm x) = 1+\mathrm{sgn}(x_1 x_2-0.25) . 
\] 
where $\mathrm{sgn}(x)$ is the sign function.
We consider the problem with both a smooth right-hand side when choosing $f = f_1$ and a non-smooth one choosing $f = f_2$.

\begin{figure}
\centering
\begin{minipage}{0.49\linewidth}
\centering
\includegraphics[width=\linewidth]{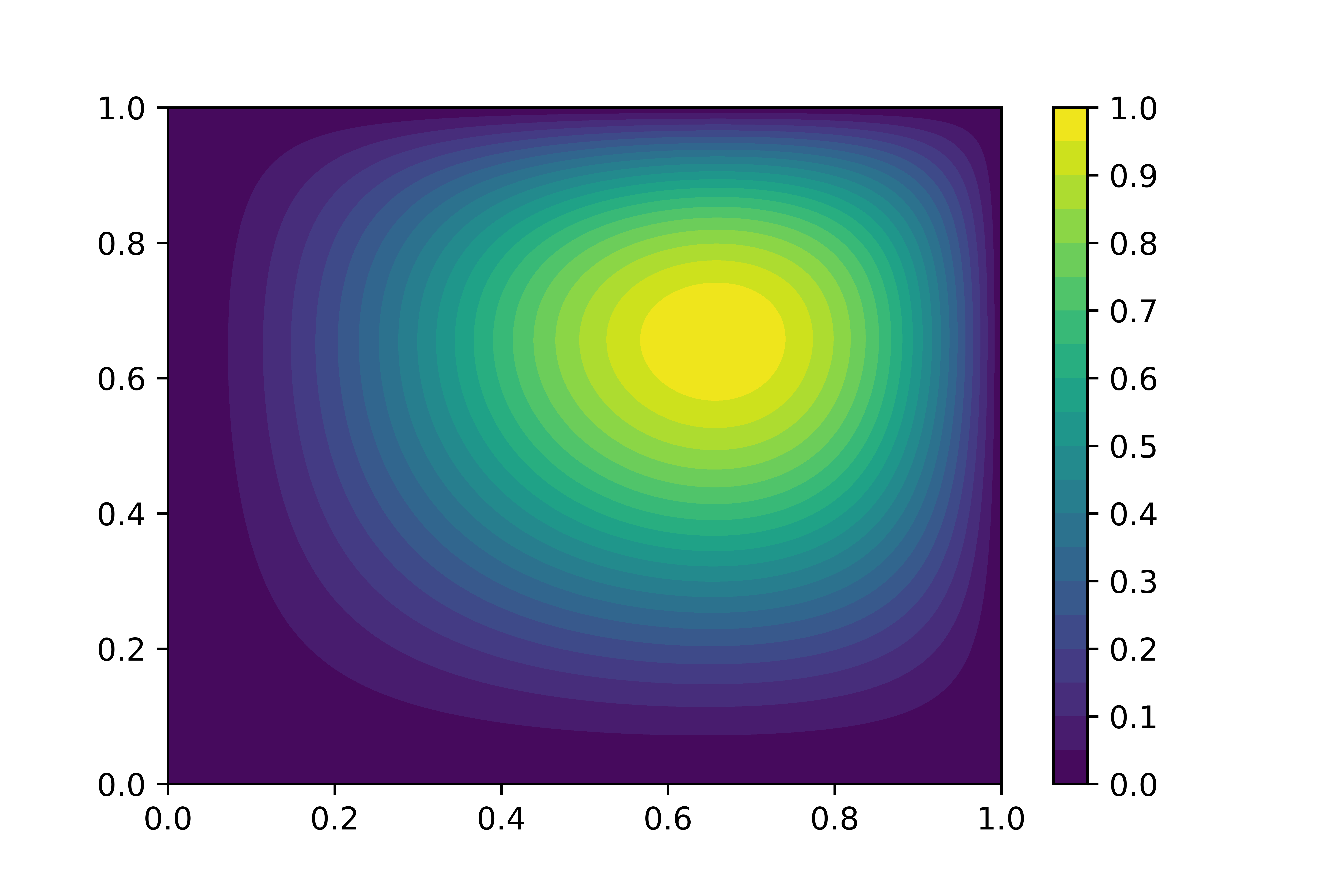}\\
$\alpha  = 0.1, \ \max u(\bm x) = 1.576239e\mathrm{-}02 $ \\
\includegraphics[width=\linewidth]{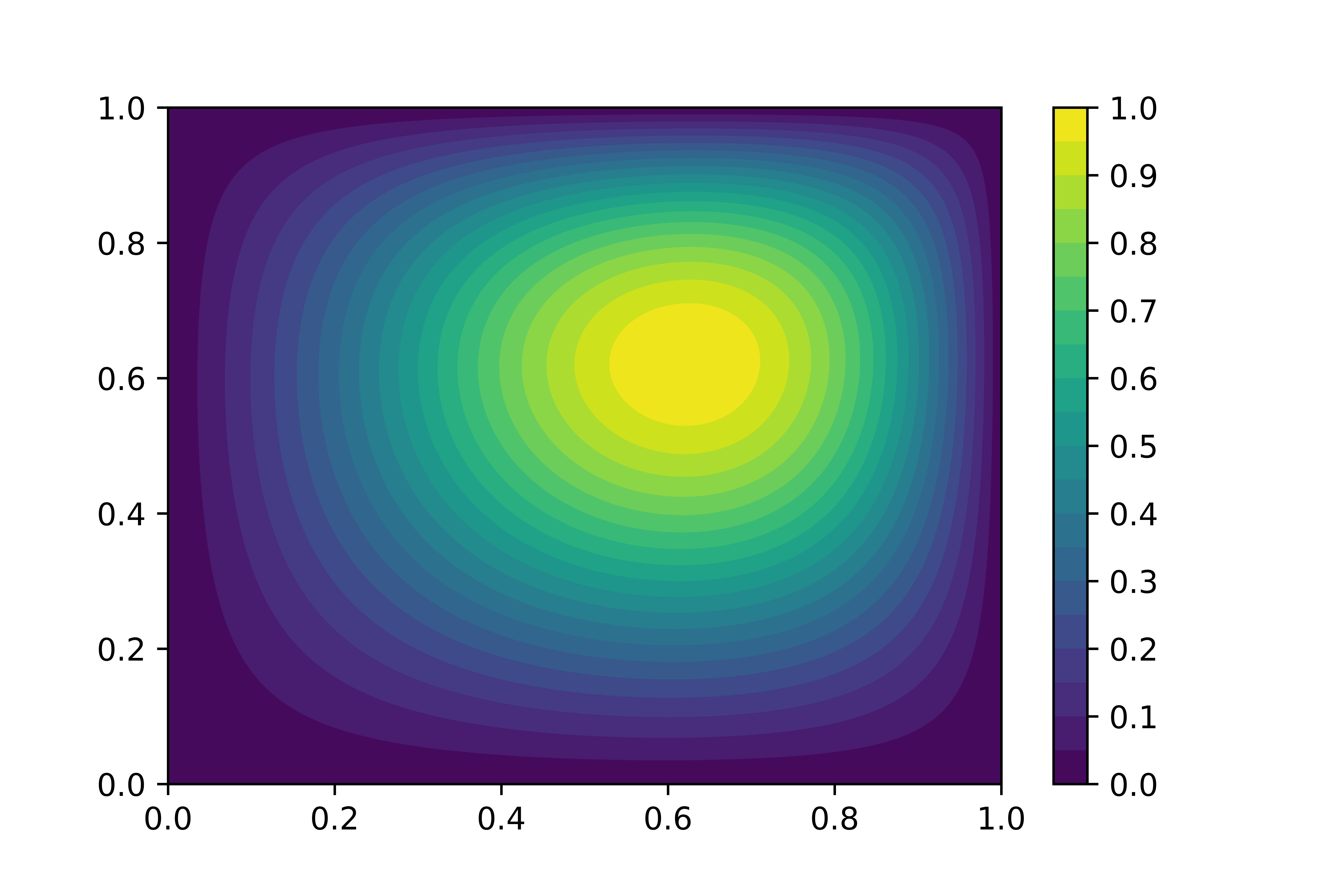}\\
$\alpha  = 0.5, \ \max u(\bm x) = 4.311910e\mathrm{-}03$ \\
\end{minipage}
\begin{minipage}{0.49\linewidth}
\centering
\includegraphics[width=\linewidth]{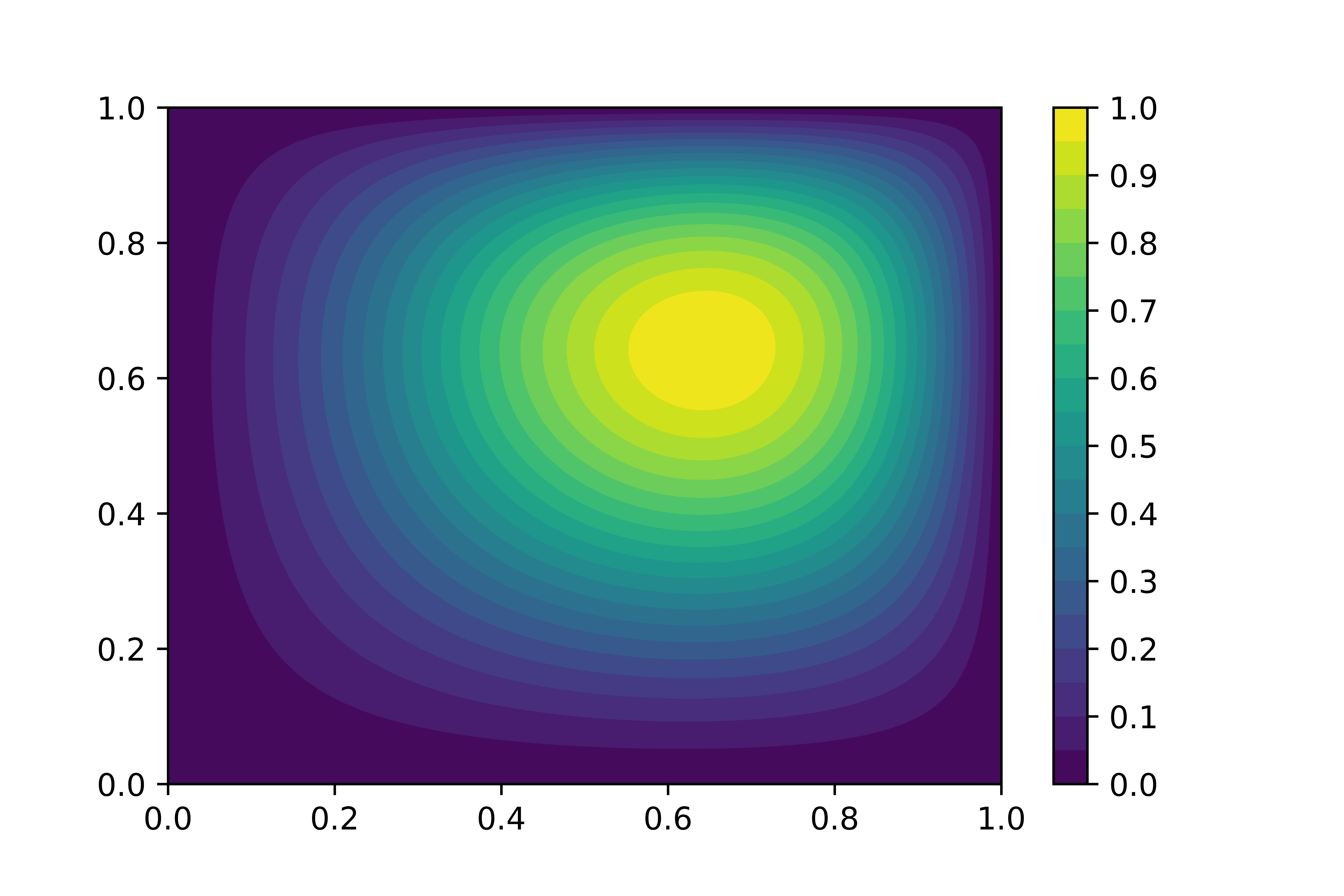}\\
$\alpha  = 0.25, \ \max u(\bm x) = 9.645557e\mathrm{-}03$ \\
\includegraphics[width=\linewidth]{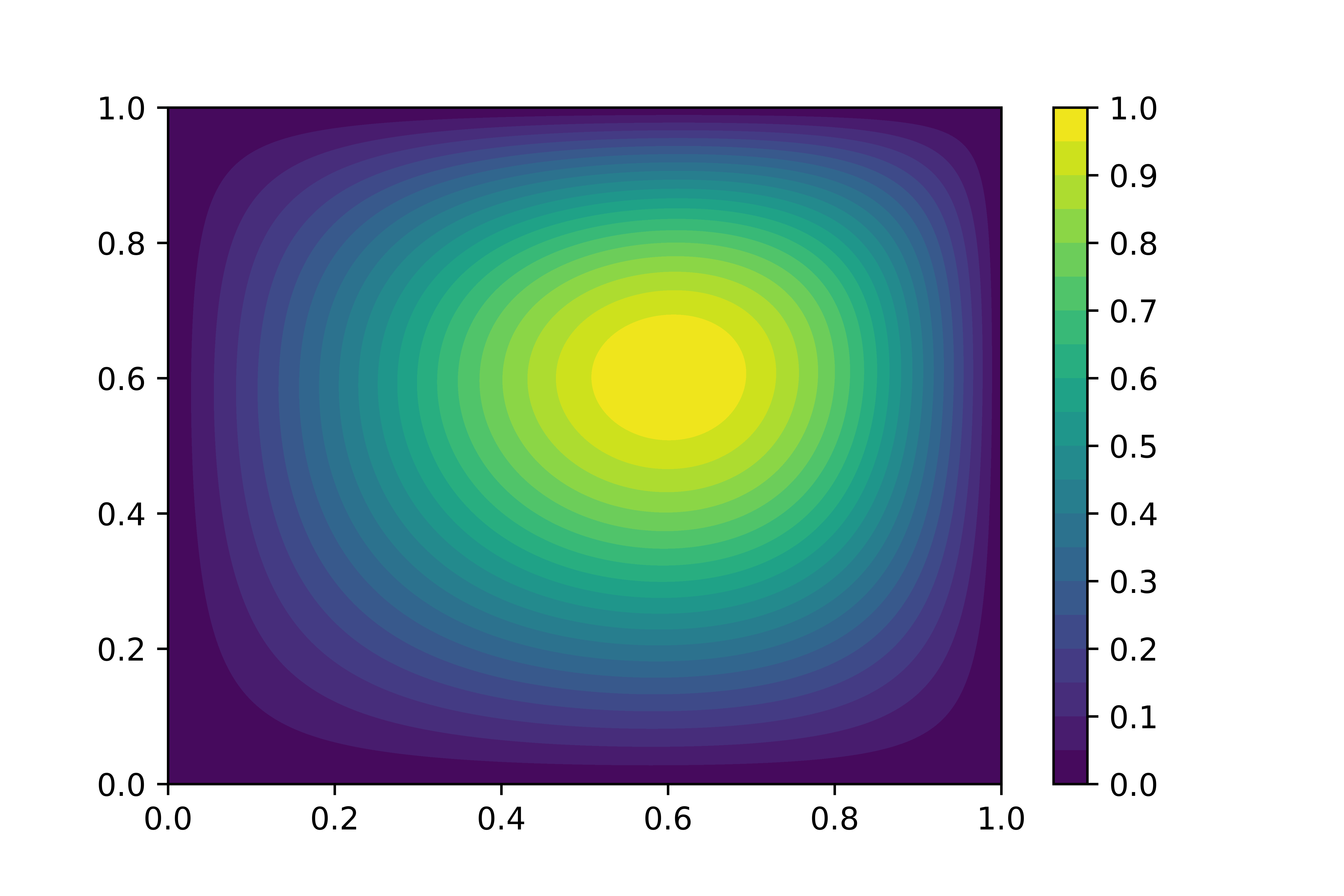}\\
$\alpha  = 0.75, \ \max u(\bm x) = 1.956855e\mathrm{-}03$ \\
\end{minipage}
\caption{Normalized solution of the problem (\ref{6}) when $f(\bm x) = f_1 (\bm x)$ on the grid $N_1 = N_2 = 256$.}
\label{f-1}
\end{figure}

\begin{figure}
\centering
\begin{minipage}{0.49\linewidth}
\centering
\includegraphics[width=\linewidth]{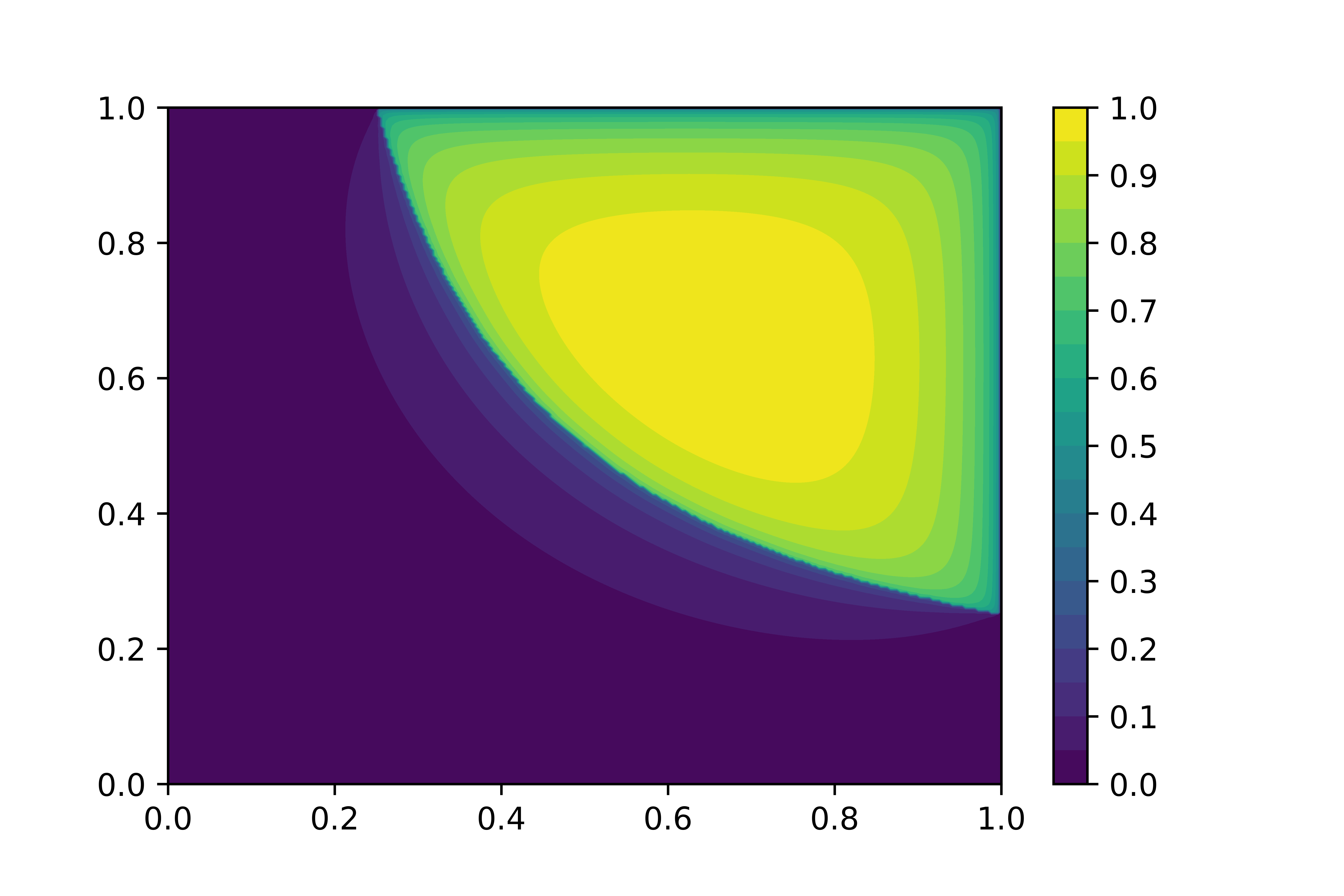}\\
$\alpha  = 0.1, \ \max u(\bm x) = 1.483479e\mathrm{+}00$ \\
\includegraphics[width=\linewidth]{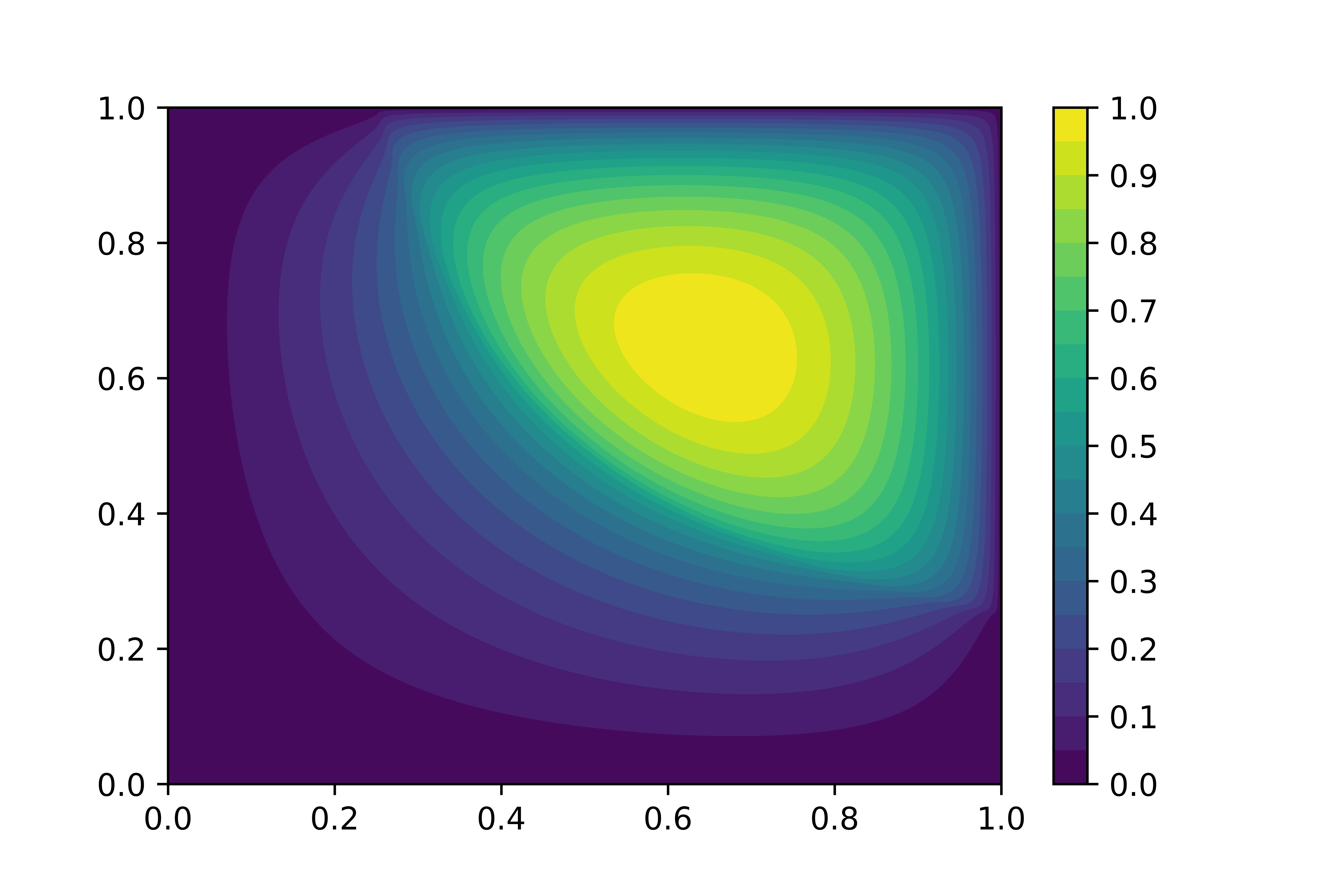}\\
$\alpha  = 0.5, \ \max u(\bm x) = 4.247955e\mathrm{-}01$ \\
\end{minipage}
\begin{minipage}{0.49\linewidth}
\centering
\includegraphics[width=\linewidth]{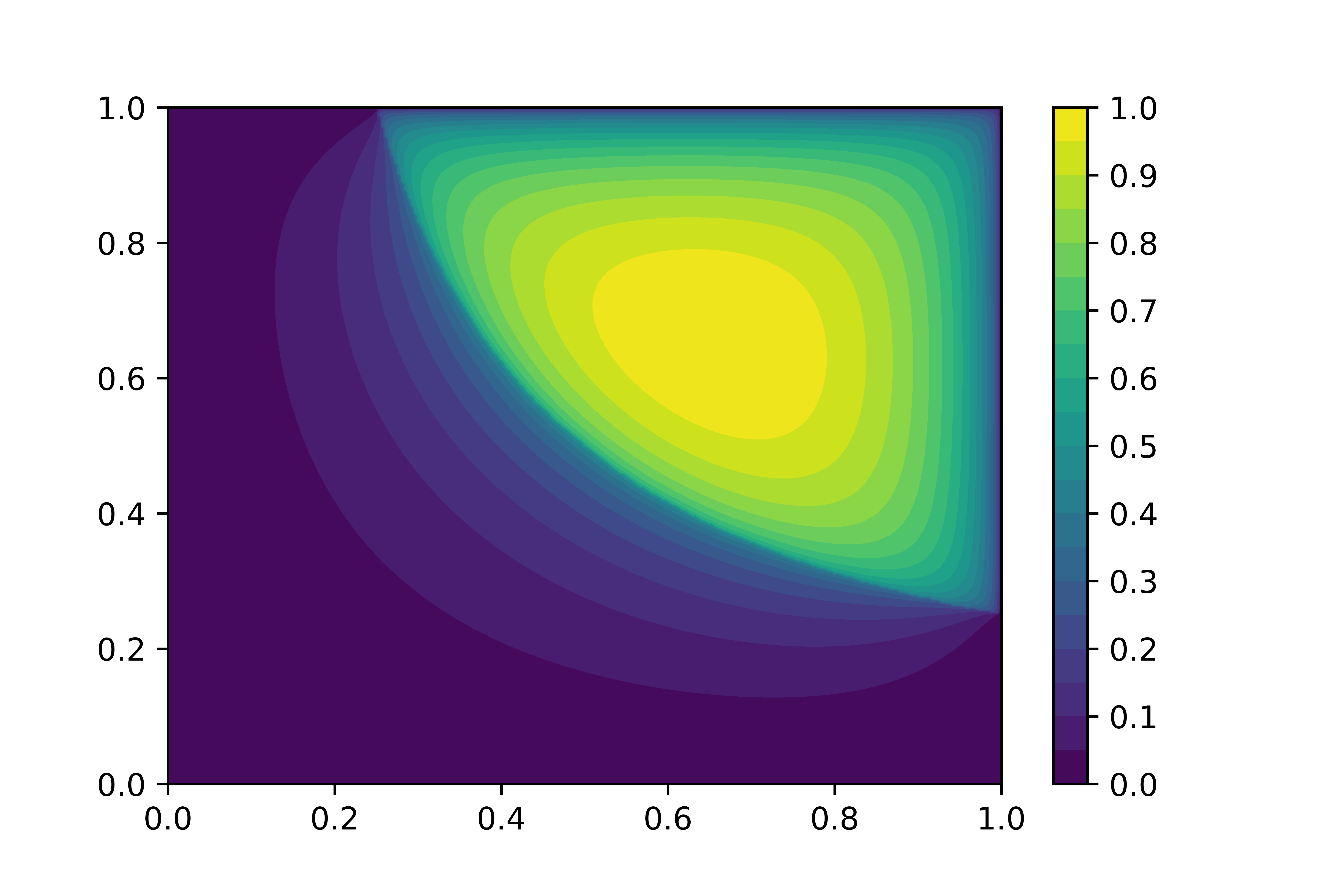}\\
$\alpha  = 0.25, \ \max u(\bm x) = 9.355516e\mathrm{-}01$ \\
\includegraphics[width=\linewidth]{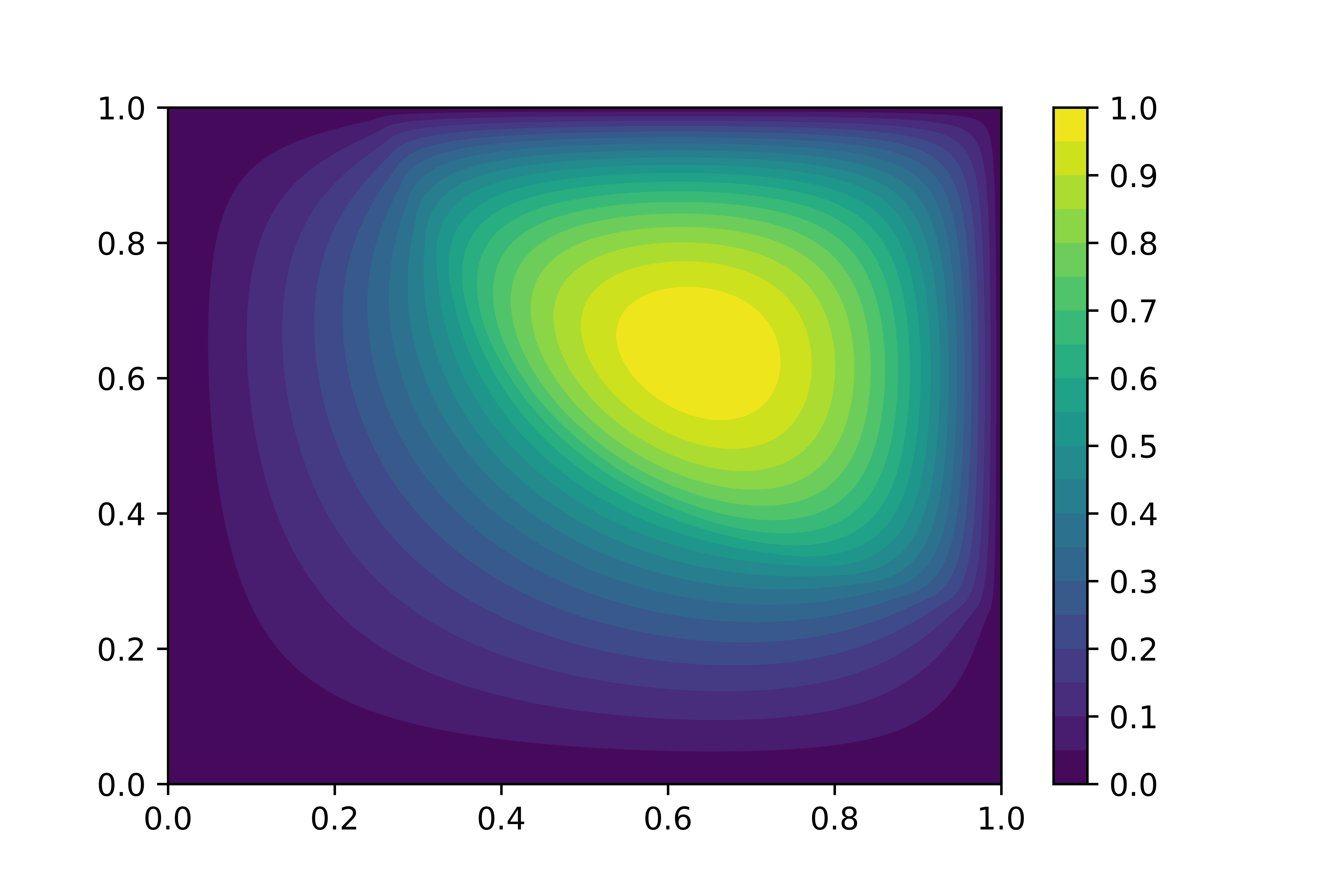}\\
$\alpha  = 0.75, \ \max u(\bm x) = 1.906235e\mathrm{-}01$ \\
\end{minipage}
\caption{Normalized solution of the problem (\ref{6}) when $f(\bm x) = f_2 (\bm x)$ on the grid $N_1 = N_2 = 256$.}
\label{f-2}
\end{figure}

The numerical solution of the problem (\ref{7}) with the right-hand side 
$b(\bm x) = f(\bm x) = f_\beta (\bm x), \ \beta = 1,2, \ \bm x \in \omega$
is depicted in Figure~\ref{f-1} and Figure~\ref{f-2} for different values of $\alpha$.
For convenience of the comparison, the following function is shown
\[
 y(\bm x) = \frac{1}{\max u(\bm x)} u(\bm x) ,
 \quad \bm x \in \omega . 
\]
For the discontinuous right-hand side ($f = f_2$, Figure~\ref{f-2}),
we observe the formation of internal and boundary layers with decreasing $\alpha$.

The relative error the approximate solution of the problem with the fractional power of the elliptic operator
is determined as follows:
\[
 \varepsilon_2 = \frac{\|\widetilde{u} -u\|}{\|u\|} ,
 \quad \varepsilon_\infty  = \frac{\|\widetilde{u}-u\|_\infty }{\|u\|_\infty } ,
 \quad \|u\|_\infty  = \max_{\bm x \in \omega} |u(\bm x)| ,
\] 
where $u$ is the exact solution of the problem (\ref{7}) and $\widetilde{u}$ is the approximate one.
The calculations were performed using the generalized Gauss-Laguerre quadrature formula with $\delta = \mu_1$.
A fairly detailed grid is used with $N_1 = N_2 = 256$.  

The dependence of the accuracy of the approximate solution of the problem (\ref{7}) on the key computational parameters
is traced according to the data in Table 2 and Table 3.
In particular, the most important is the dependence on the number of nodes of the quadrature formula $M$.
Note that as the parameter $p$ increases, the accuracy of the approximate solution does not decrease.
Therefore, we choose $p = 0$. Comparison of Table 2 and Table 3 shows a significant drop in accuracy
for a problem with a non-smooth right side.    

\begin{center}
\begin{table}[htp]
\label{t-2}
\caption{Solution error for the problem (\ref{7}) with $f(\bm x) = f_1 (\bm x)$ at various parameter values $p$}
\centering
\begin{tabular}{cccccccc}
\hline
  $p$  & $M$   &    error  &   $\alpha = 0.1$         &     $\alpha = 0.25$   &    $\alpha = 0.5$   &     $\alpha = 0.75$      &    $\alpha = 0.9$   \\
\hline
0  &     25   &      $\varepsilon_2$          &  1.867333e-07     &  2.234522e-07     &    9.979937e-08    &  3.388432e-08    &  1.713416e-08  \\
   &          &      $\varepsilon_{\infty}$   &  1.101566e-03     &  1.066598e-03     &    6.414770e-04    &  3.516336e-04    &  2.428239e-04  \\
   &     50   &      $\varepsilon_2$          &  2.826563e-08     &  2.756374e-08     &    8.744432e-09    &  2.113860e-09    &  8.727410e-10  \\
   &          &      $\varepsilon_{\infty}$   &  5.127673e-04     &  4.435298e-04     &    2.247211e-04    &  1.025590e-04    &  6.344683e-05  \\
   &    100   &      $\varepsilon_2$          &  4.257879e-09     &  3.399113e-09     &    7.664039e-10    &  1.315780e-10    &  4.425344e-11  \\
   &          &      $\varepsilon_{\infty}$   &  2.358176e-04     &  1.864016e-04     &    7.926752e-05    &  3.057283e-05    &  1.697664e-05  \\
\hline   
1  &     25   &      $\varepsilon_2$          &  1.042484e-05     &  4.322339e-06     &    1.048142e-06    &  2.702530e-07    &  1.233633e-07  \\
   &          &      $\varepsilon_{\infty}$   &  5.917636e-03     &  3.726284e-03     &    1.777896e-03    &  8.772466e-04    &  5.808148e-04  \\
   &     50   &      $\varepsilon_2$          &  1.604299e-06     &  5.427969e-07     &    9.390374e-08    &  1.725236e-08    &  6.414784e-09  \\
   &          &      $\varepsilon_{\infty}$   &  2.725576e-03     &  1.538243e-03     &    6.132594e-04    &  2.543889e-04    &  1.527295e-04  \\
   &    100   &      $\varepsilon_2$          &  2.460128e-07     &  6.779333e-08     &    8.335361e-09    &  1.089103e-09    &  3.301012e-10  \\
   &          &      $\varepsilon_{\infty}$   &  1.267140e-03     &  6.489867e-04     &    2.177537e-04    &  7.583209e-05    &  4.092155e-05  \\
\hline     
2  &     25   &      $\varepsilon_2$          &  3.445381e-05     &  1.438506e-05     &    3.540165e-06    &  9.253987e-07    &  4.236814e-07  \\
   &          &      $\varepsilon_{\infty}$   &  9.668806e-03     &  6.152098e-03     &    2.969534e-03    &  1.475692e-03    &  9.797893e-04  \\
   &     50   &      $\varepsilon_2$          &  5.381860e-06     &  1.832602e-06     &    3.200697e-07    &  5.929701e-08    &  2.213631e-08  \\
   &          &      $\varepsilon_{\infty}$   &  4.393399e-03     &  2.526119e-03     &    1.037259e-03    &  4.408139e-04    &  2.672836e-04  \\
   &    100   &      $\varepsilon_2$          &  8.326382e-07     &  2.311031e-07     &    2.870886e-08    &  3.782843e-09    &  1.151544e-09  \\
   &          &      $\varepsilon_{\infty}$   &  2.059608e-03     &  1.068331e-03     &    3.650104e-04    &  1.287746e-04    &  6.980999e-05  \\
\hline
\end{tabular}
\end{table}
\end{center}

\begin{center}
\begin{table}[htp]
\label{t-3}
\caption{Solution error for the problem (\ref{7}) with $f(\bm x) = f_2 (\bm x)$ at various parameter values $p$}
\centering
\begin{tabular}{cccccccc}
\hline
  $p$  & $M$   &    error  &   $\alpha = 0.1$         &     $\alpha = 0.25$   &    $\alpha = 0.5$   &     $\alpha = 0.75$      &    $\alpha = 0.9$   \\
\hline
0  &     25   &      $\varepsilon_2$          &  1.084357e-02     &  3.938291e-03     &    6.485652e-04    &  1.177557e-04    &  4.442441e-05  \\
   &          &      $\varepsilon_{\infty}$   &  4.370840e-01     &  1.905319e-01     &    6.248551e-02    &  2.430478e-02    &  1.442032e-02  \\
   &     50   &      $\varepsilon_2$          &  6.312810e-03     &  1.936705e-03     &    2.296383e-04    &  2.977748e-05    &  9.176919e-06  \\
   &          &      $\varepsilon_{\infty}$   &  4.005614e-01     &  1.588500e-01     &    4.413926e-02    &  1.444312e-02    &  7.761052e-03  \\
   &    100   &      $\varepsilon_2$          &  3.552338e-03     &  9.377168e-04     &    8.059058e-05    &  7.455089e-06    &  1.873778e-06  \\
   &          &      $\varepsilon_{\infty}$   &  3.725840e-01     &  1.330906e-01     &    3.095324e-02    &  8.540190e-03    &  4.166878e-03  \\
\hline   
1  &     25   &      $\varepsilon_2$          &  4.076650e-02     &  1.271564e-02     &    2.054384e-03    &  3.768885e-04    &  1.428314e-04  \\
   &          &      $\varepsilon_{\infty}$   &  5.324292e-01     &  2.576112e-01     &    9.243630e-02    &  3.763498e-02    &  2.275878e-02  \\
   &     50   &      $\varepsilon_2$          &  2.485088e-02     &  6.386131e-03     &    7.389697e-04    &  9.694268e-05    &  3.010355e-05  \\
   &          &      $\varepsilon_{\infty}$   &  4.969681e-01     &  2.169114e-01     &    6.579678e-02    &  2.259551e-02    &  1.239358e-02  \\
   &    100   &      $\varepsilon_2$          &  1.494239e-02     &  3.170517e-03     &    2.626296e-04    &  2.456552e-05    &  6.228987e-06  \\
   &          &      $\varepsilon_{\infty}$   &  4.616624e-01     &  1.812017e-01     &    4.657524e-02    &  1.352842e-02    &  6.689115e-03  \\   
\hline     
2  &     25   &      $\varepsilon_2$          &  5.766283e-02     &  1.970018e-02     &    3.575837e-03    &  7.088258e-04    &  2.777896e-04  \\
   &          &      $\varepsilon_{\infty}$   &  5.589886e-01     &  2.869214e-01     &    1.107642e-01    &  4.739887e-02    &  2.933699e-02  \\
   &     50   &      $\varepsilon_2$          &  3.549610e-02     &  9.985667e-03     &    1.305410e-03    &  1.877383e-04    &  6.099850e-05  \\
   &          &      $\varepsilon_{\infty}$   &  5.225493e-01     &  2.425210e-01     &    7.941726e-02    &  2.898929e-02    &  1.636029e-02  \\
   &    100   &      $\varepsilon_2$          &  2.157698e-02     &  5.001624e-03     &    4.685264e-04    &  4.819297e-05    &  1.280421e-05  \\
   &          &      $\varepsilon_{\infty}$   &  4.867483e-01     &  2.038174e-01     &    5.642903e-02    &  1.737444e-02    &  8.854464e-03  \\
\hline
\end{tabular}
\end{table}
\end{center}

The dependence of the accuracy of the approximate solutions on the discrete problem (\ref{7}) of the number of nodes in space is of great interest as well.
The calculation results for various values of $N = N_1 = N_2$ are
presented in Table 4 and Table 5.
The dimension of a discrete problem practically does not affect the accuracy of
quadrature formula for an approximate solution of fractional elliptical BVP.
  
\begin{center}
\begin{table}[htp]
\label{t-4}
\caption{Solution error for the problem (\ref{7}) with $f(\bm x) = f_1 (\bm x)$ at various parameter values $N = N_1 = N_2$}
\centering
\begin{tabular}{cccccccc}
\hline
  $N$  & $M$   &    error  &   $\alpha = 0.1$         &     $\alpha = 0.25$   &    $\alpha = 0.5$   &     $\alpha = 0.75$      &    $\alpha = 0.9$   \\
\hline
32 &     25   &      $\varepsilon_2$          &  8.409942e-08     &  1.456405e-07     &    7.988626e-08    &  2.925977e-08    &  1.517098e-08  \\
   &          &      $\varepsilon_{\infty}$   &  6.600418e-04     &  8.573971e-04     &    5.919345e-04    &  3.250000e-04    &  2.188000e-04  \\
   &     50   &      $\varepsilon_2$          &  5.068313e-09     &  9.355974e-09     &    4.825193e-09    &  1.451950e-09    &  6.457840e-10  \\
   &          &      $\varepsilon_{\infty}$   &  1.553199e-04     &  2.178892e-04     &    1.602107e-04    &  8.719835e-05    &  5.690694e-05  \\
   &    100   &      $\varepsilon_2$          &  1.238361e-10     &  2.538166e-10     &    1.438117e-10    &  4.265139e-11    &  1.791286e-11  \\
   &          &      $\varepsilon_{\infty}$   &  2.100015e-05     &  3.149206e-05     &    2.531976e-05    &  1.457506e-05    &  9.665245e-06  \\
\hline   
64 &     25   &      $\varepsilon_2$          &  1.610774e-07     &  2.073936e-07     &    9.562513e-08    &  3.287397e-08    &  1.669788e-08  \\
   &          &      $\varepsilon_{\infty}$   &  1.032142e-03     &  9.701921e-04     &    6.310096e-04    &  3.372346e-04    &  2.391175e-04  \\
   &     50   &      $\varepsilon_2$          &  2.002888e-08     &  2.324099e-08     &    7.969464e-09    &  1.982252e-09    &  8.265304e-10  \\
   &          &      $\varepsilon_{\infty}$   &  4.305003e-04     &  4.156202e-04     &    1.998006e-04    &  9.833823e-05    &  6.224204e-05  \\
   &    100   &      $\varepsilon_2$          &  1.943737e-09     &  2.230604e-09     &    6.161669e-10    &  1.139644e-10    &  3.927857e-11  \\
   &          &      $\varepsilon_{\infty}$   &  1.420510e-04     &  1.499068e-04     &    7.332841e-05    &  2.840727e-05    &  1.554351e-05  \\
\hline     
128&     25   &      $\varepsilon_2$          &  1.820715e-07     &  2.204179e-07     &    9.898763e-08    &  3.368579e-08    &  1.704815e-08  \\
   &          &      $\varepsilon_{\infty}$   &  1.087668e-03     &  1.060296e-03     &    6.395206e-04    &  3.508900e-04    &  2.412719e-04  \\
   &     50   &      $\varepsilon_2$          &  2.676261e-08     &  2.678691e-08     &    8.598566e-09    &  2.088511e-09    &  8.637732e-10  \\
   &          &      $\varepsilon_{\infty}$   &  4.754600e-04     &  4.383331e-04     &    2.231665e-04    &  1.006521e-04    &  6.322556e-05  \\
   &    100   &      $\varepsilon_2$          &  3.761240e-09     &  3.195480e-09     &    7.399354e-10    &  1.283428e-10    &  4.332366e-11  \\
   &          &      $\varepsilon_{\infty}$   &  2.236805e-04     &  1.701325e-04     &    7.822322e-05    &  2.923203e-05    &  1.683085e-05  \\ 
\hline
\end{tabular}
\end{table}
\end{center}

\begin{center}
\begin{table}[htp]
\label{t-5}
\caption{Solution error for the problem (\ref{7}) with $f(\bm x) = f_2 (\bm x)$ at various parameter values $N = N_1 = N_2$}
\centering
\begin{tabular}{cccccccc}
\hline
  $N$  & $M$   &    error  &   $\alpha = 0.1$         &     $\alpha = 0.25$   &    $\alpha = 0.5$   &     $\alpha = 0.75$      &    $\alpha = 0.9$   \\
\hline
32 &     25   &      $\varepsilon_2$          &  2.006090e-03     &  1.828621e-03     &    4.691900e-04    &  9.911544e-05    &  3.910657e-05  \\
   &          &      $\varepsilon_{\infty}$   &  1.187981e-01     &  1.098096e-01     &    5.279046e-02    &  2.143857e-02    &  1.226111e-02  \\
   &     50   &      $\varepsilon_2$          &  3.760510e-04     &  4.506979e-04     &    1.189540e-04    &  2.059487e-05    &  6.933095e-06  \\
   &          &      $\varepsilon_{\infty}$   &  5.566576e-02     &  5.283872e-02     &    2.524998e-02    &  1.103353e-02    &  6.396780e-03  \\
   &    100   &      $\varepsilon_2$          &  3.792739e-05     &  5.832591e-05     &    1.907976e-05    &  3.269577e-06    &  1.009050e-06  \\
   &          &      $\varepsilon_{\infty}$   &  1.965240e-02     &  2.111303e-02     &    1.054680e-02    &  4.130081e-03    &  2.183357e-03  \\
\hline   
64 &     25   &      $\varepsilon_2$          &  6.360886e-03     &  3.165742e-03     &    5.971673e-04    &  1.123797e-04    &  4.279690e-05  \\
   &          &      $\varepsilon_{\infty}$   &  3.258161e-01     &  1.792333e-01     &    5.846466e-02    &  2.324461e-02    &  1.372679e-02  \\
   &     50   &      $\varepsilon_2$          &  2.574703e-03     &  1.312433e-03     &    1.966289e-04    &  2.738334e-05    &  8.597362e-06  \\
   &          &      $\varepsilon_{\infty}$   &  2.021267e-01     &  1.333749e-01     &    3.999555e-02    &  1.290906e-02    &  7.372346e-03  \\
   &    100   &      $\varepsilon_2$          &  7.619762e-04     &  4.606884e-04     &    5.933321e-05    &  6.302869e-06    &  1.650515e-06  \\
   &          &      $\varepsilon_{\infty}$   &  1.036411e-01     &  7.774238e-02     &    2.644230e-02    &  7.584779e-03    &  3.532676e-03  \\
\hline     
128 &     25   &      $\varepsilon_2$          &  9.378730e-03     &  3.740118e-03     &    6.383832e-04    &  1.168124e-04    &  4.414578e-05  \\
   &          &      $\varepsilon_{\infty}$   &  4.279399e-01     &  1.880510e-01     &    6.172671e-02    &  2.397168e-02    &  1.435631e-02  \\
   &     50   &      $\varepsilon_2$          &  4.994286e-03     &  1.766193e-03     &    2.226535e-04    &  2.933208e-05    &  9.074022e-06  \\
   &          &      $\varepsilon_{\infty}$   &  3.718976e-01     &  1.502341e-01     &    4.274832e-02    &  1.427241e-02    &  7.658653e-03  \\
   &    100   &      $\varepsilon_2$          &  2.419279e-03     &  7.963929e-04     &    7.576465e-05    &  7.224747e-06    &  1.831256e-06  \\
   &          &      $\varepsilon_{\infty}$   &  2.838417e-01     &  1.268924e-01     &    2.938177e-02    &  8.281267e-03    &  4.000634e-03  \\
\hline
\end{tabular}
\end{table}
\end{center}

\section{Conclusions}\label{sec:6} 

\begin{enumerate}
 \item To solve the problem with a fractional power of a positive definite operator, we use the known integral representation of the solution by solving the Cauchy problem. The possibilities of such representation are expanded due to the transition to the problem with the new initial conditions. The key dependence of the solution on the time of the auxiliary non-stationary problem, associated with the constant positive definiteness operator, is defined for constructing optimal Gaussian quadratures.
 \item An approximate solution is based on the use of the generalized Gauss-Laguerre quadrature formula. The accuracy of the quadrature formula is studied on the representative parametric integration problem depending on the key parameters: degrees $0 <\alpha <1$ and the number of nodes of the quadrature formula.
 \item The model problem with the fractional power of the Laplace operator with a smooth and discontinuous right-hand side is considered in a rectangle. Its approximate solution is obtained based on the integral representation by solving the parabolic Cauchy problems. The results of the numerical experiments on the influence of the smoothness of the right-hand side on the accuracy of the approximate solutions using a different number of integration nodes are presented.
\end{enumerate}


\begin{thebibliography}{10}
\expandafter\ifx\csname url\endcsname\relax
  \def\url#1{\texttt{#1}}\fi
\expandafter\ifx\csname urlprefix\endcsname\relax\def\urlprefix{URL }\fi
\expandafter\ifx\csname href\endcsname\relax
  \def\href#1#2{#2} \def\path#1{#1}\fi

\bibitem{baleanu2012fractional}
D.~Baleanu, Fractional Calculus: Models and Numerical Methods, World
  Scientific, New York, 2012.

\bibitem{uchaikin}
V.~V. Uchaikin, Fractional Derivatives for Physicists and Engineers, Higher
  Education Press, 2013.

\bibitem{Pozrikidis16}
C.~Pozrikidis, The Fractional Laplacian, CRC Press, 2018.

\bibitem{higham2008functions}
N.~J. Higham, Functions of Matrices: Theory and Computation, SIAM,
  Philadelphia, 2008.

\bibitem{bonito2017}
A.~Bonito, J.~P. Borthagaray, R.~H. Nochetto, E.~Otarola, A.~J. Salgado,
  Numerical methods for fractional diffusion, Computing and Visualization in
  Science 19~(5-6) (2018) 19--46.

\bibitem{stahl2003best}
H.~R. Stahl, Best uniform rational approximation of $x^\alpha$ on [0, 1], Acta
  Mathematica 190~(2) (2003) 241--306.

\bibitem{harizanov2018optimal}
S.~Harizanov, R.~Lazarov, S.~Margenov, P.~Marinov, Y.~Vutov, Optimal solvers
  for linear systems with fractional powers of sparse spd matrices, Numerical
  Linear Algebra with Applications 25~(5) (2018) e2167.

\bibitem{HarizanovArxive2019}
S.~Harizanov, R.~Lazarov, P.~Marinov, S.~Margenov, J.~Pasciak, Analysis of
  numerical methods for spectral fractional elliptic equations based on the
  best uniform rational approximation, arXiv:1905.08155.

\bibitem{balakrishnan1960fractional}
A.~V. Balakrishnan, Fractional powers of closed operators and the semigroups
  generated by them., Pacific Journal of Mathematics 10~(2) (1960) 419--437.

\bibitem{birman1987spectral}
M.~S. Birman, M.~Z. Solomjak, Spectral Theory of Self-adjoint Operators in
  Hilbert Space, Kluwer academic publishers, 1987.

\bibitem{bonito2015numerical}
A.~Bonito, J.~Pasciak, Numerical approximation of fractional powers of elliptic
  operators, Mathematics of Computation 84~(295) (2015) 2083--2110.

\bibitem{AcetoNovat}
L.~Aceto, P.~Novati, Rational approximation to the fractional {L}aplacian
  operator in reaction-diffusion problems, SIAM Journal on Scientific Computing
  39~(1) (2017) A214--A228.

\bibitem{Aceto2019}
L.~Aceto, P.~Novati, Rational approximations to fractional powers of
  self-adjoint positive operators authors, Numerische Mathematik 143~(1) (2019)
  1--16.

\bibitem{vab2019-10838}
P.~N. Vabishchevich, Approximation of a fractional power of an elliptic
  operator, arXiv:1905.10838.

\bibitem{Caffarelli}
L.~Caffarelli, L.~Silvestre, An extension problem related to the fractional
  {L}aplacian, Communications in Partial Differential Equations 32~(8) (2007)
  1245--1260.

\bibitem{nochetto2015pde}
R.~H. Nochetto, E.~Ot{\'a}rola, A.~J. Salgado, A {PDE} approach to fractional
  diffusion in general domains: a priori error analysis, Foundations of
  Computational Mathematics 15~(3) (2015) 733--791.

\bibitem{vabishchevich2014numerical}
P.~N. Vabishchevich, Numerically solving an equation for fractional powers of
  elliptic operators, Journal of Computational Physics 282~(1) (2015) 289--302.

\bibitem{duan2018numerical}
B.~Duan, R.~Lazarov, J.~Pasciak, Numerical approximation of fractional powers
  of elliptic operators, IMA Journal of Numerical Analysis~(drz013) (2019)
  1--26.
\newblock \href {http://dx.doi.org/10.1093/imanum/drz013}
  {\path{doi:10.1093/imanum/drz013}}.

\bibitem{ciegisvab2019-00201}
R.~Ciegis, P.~Vabishchevich, High order numerical schemes for solving
  fractional powers of elliptic operators, Computers \& Mathematics with
  Applications (2019) 1--11.~In Press.
\newblock \href {http://dx.doi.org/10.1016/j.camwa.2019.08.012}
  {\path{doi:10.1016/j.camwa.2019.08.012}}.

\bibitem{stinga2019user}
P.~R. Stinga, User’s guide to the fractional {L}aplacian and the method of
  semigroups, in: Handbook of Fractional Calculus and Applications. Volume 2:
  Fractional Differential Equations, de Gruyter, 2019, pp. 235--266.

\bibitem{cusimano2018discretizations}
N.~Cusimano, F.~del Teso, L.~Gerardo-Giorda, G.~Pagnini, Discretizations of the
  spectral fractional {L}aplacian on general domains with {D}irichlet,
  {N}eumann, and {R}obin boundary conditions, SIAM Journal on Numerical
  Analysis 56~(3) (2018) 1243--1272.

\bibitem{cusimano2018numerical}
N.~Cusimano, F.~del Teso, L.~Gerardo-Giorda, Numerical approximations for
  fractional elliptic equations via the method of semigroups, arXiv:1812.01518.

\bibitem{davis2007methods}
P.~J. Davis, P.~Rabinowitz, Methods of Numerical Integration, Dover
  Publications, 2007.

\bibitem{Samarskii1989}
A.~A. Samarskii, The Theory of Difference Schemes, Marcel Dekker, New York,
  2001.

\bibitem{SamarskiiNikolaev1978}
A.~A. Samarskii, E.~S. Nikolaev, Numerical Methods for Grid Equations. Vol. I,
  II, Birkhauser Verlag, Basel, 1989.

\end{thebibliography}
\end{document}